\documentclass
{amsart}
\usepackage{amsmath, amsthm, amssymb,mathrsfs,amsfonts,}
\usepackage{mathtools}
\usepackage{relsize}
\usepackage{enumitem}
\usepackage{pgffor}

\usepackage[
]
{hyperref}

\usepackage{graphicx}

\newcommand\scirc[1][.45]{\mathbin{\vcenter{\hbox{\scalebox{#1}{$\circ$}}}}}

\usepackage{cases} 
\usepackage{tikz}

\newcommand{\Be}{\begin{equation}}
\newcommand{\Ee}{\end{equation}}
\newcommand{\Bea}{\begin{eqnarray}}
\newcommand{\Eea}{\end{eqnarray}}
\newcommand{\Bel}{\begin{align}}
\newcommand{\Eel}{\end{align}}
\newcommand{\Beas}{\begin{eqnarray*}}
	\newcommand{\Eeas}{\end{eqnarray*}}
\newcommand{\Benu}{\begin{enumerate}}
	\newcommand{\Eenu}{\end{enumerate}}
\newcommand{\Bi}{\begin{itemize}}
	\newcommand{\Ei}{\end{itemize}}

\numberwithin{equation}{section}
\newcommand{\dist} {\text{dist\! }}

\theoremstyle{plain}
\newtheorem{thm}{Theorem}[section]
\newtheorem{cor}[thm]{Corollary}

\newtheorem{lem}[thm]{Lemma}
\newtheorem{prop}[thm]{Proposition}

\newtheorem{conj}[thm]{Conjecture}
  
\newtheorem{defn}[thm]{Definition}

\definecolor{ao}{rgb}{0, 0.5, 0}

\newcommand{\delc}{\delta_{\!\, \mathsmaller{\mathsmaller {\scirc}}}}
\newcommand{\sz}{s_{\!\,\mathsmaller{\mathsmaller{\scirc}}}}

\newcommand{\fa}{\mathfrak a}

\DeclareMathOperator{\Vol}{Vol}
\DeclareMathOperator{\sspan}{span}
\DeclareMathOperator{\supp}{supp}

\begin{document}

\title[Regularity of restricted X-ray transforms]{
Sharp Sobolev regularity  of  restricted \\  X-ray transforms
}

\author[H. Ko]{Hyerim Ko}
\author[S. Lee]{Sanghyuk Lee}
\author[S. Oh]{Sewook  Oh}

\address{Department of Mathematical Sciences and RIM, Seoul National University, Seoul 08826, Republic of Korea}
\email{kohr@snu.ac.kr}
\email{shklee@snu.ac.kr}

\address{School of Mathematics, Korea Institute for Advanced Study, 85 Hoegiro Dongdaemun-gu, Seoul 02455, Republic of Korea}
\email{sewookoh@kias.re.kr}


\begin{abstract}    
We study $L^p$-Sobolev regularity estimate for  the  restricted  X-ray transforms generated by  nondegenerate curves.
 Making use of  the inductive strategy in the recent work by the authors \cite{KLO}, we establish the sharp  $L^p$-regularity estimates for the  restricted X-ray transforms in $\mathbb R^{d+1}$, $d\ge 3$. 
 This extends  the result due to Pramanik and Seeger \cite{PS06} in $\mathbb R^3$ to every dimension.
\end{abstract}

\maketitle 

\section{Introduction}
Let $\gamma$ be a smooth curve  from $I=[-1,1]$ to $\mathbb R^{d}$.
We consider 
\[
\mathfrak R f(x,s)=\psi(s)\int f(x+t\gamma(s),t)\chi(t)dt,  \quad f\in \mathcal S(\mathbb{R}^{d+1}),
\]
where $\psi$ and $\chi$ are smooth functions supported in the interiors of the intervals $I$ and  $[1,2]$, respectively. The operator 
$\mathfrak Rf$ is referred to as the restriction of X-ray transform to the line complex  generated in the direction $(\gamma(s),1)$. 
We say  $\gamma$ is  nondegenerate  if
\begin{align}\label{nonv}
\det( \gamma'(s), \dots, \gamma^{(d)}(s) )\neq0,
\quad  \forall s\in I.
\end{align}

The operator $\mathfrak Rf$ is a model case of the general class of restricted X-ray transforms 
(see \cite{GG, GS00, GU1, GU2, GU3}). 
Especially in $\mathbb R^3$, under the nondegeneracy assumption \eqref{nonv}, $\mathfrak Rf$  is 
a typical example of Fourier integral operators with one-sided fold singularity  (\cite{GS94}). 
Regularity properties of $\mathfrak Rf$ have been studied in terms of $L^p$ improving and $L^p$ Sobolev regularity estimates. 
$L^p$ improving property of $\mathfrak R$  is well understood by now (\cite{GS98, GSW, Ob, La}). The problem was, in fact,   considered in a more general framework: 
$L^p$--$L_s^q(L_x^r)$ estimates for $\mathfrak R$ were studied by some authors (see, for example,  \cite{W, E, CE1}) and the estimates on the optimal range of $p,q$ were established 
 except for some endpoint cases.   (See also \cite{CE2, TW, G,  DS1, DS2} for related results.)

The $L^2$--$L_{1/(2d)}^2$ bound on $\mathfrak R$ is easy to obtain via $TT^*$ argument and van der Corput's lemma (\cite{GS00}) 
(see also \cite{GU3, GS94} for the sharp $L^2$ Sobolev estimates for general class of operators). 
Interpolation  between this  and the trivial $L^\infty$ estimate shows that $\mathfrak R$ is bounded from  $L^p$ to $L_{1/(pd)}^p$  for $p\ge 2$. This is  optimal in that $L^p$--$L_{\alpha}^p$ estimate fails if $\alpha > 1/(pd)$  (see Proposition \ref{nece} below). However, when $p< 2$,  the sharp $L^p$ regularity estimate  is less straightforward.   
Such estimate was not known until recently. When $d=2$,  the optimal $L^p-L_{1/p'}^p$ estimate was established for 
$1< p<4/3$ by Pramanik and Seeger's conditional result \cite{PS06} 
and the sharp decoupling inequality for the cone $\subset \mathbb R^3$ due to 
Bourgain and Demeter \cite{BD}.  Those estimates and interpolation  give the sharp $L^p$--$L_{1/(2d+)}^p$ estimate  for $4/3\le p<2$ but 
 the endpoint 
$L^p$--$L_{1/(2d)}^p$ estimate remains open. (See Conjecture \ref{sobolev} below.) In $\mathbb R^3$ the result has been  extended to more general operators. In fact, 
Pramanik and Seeger \cite{PS20} obtained the sharp $L^p$ regularity estimates for  Fourier integral operator with folding canonical relation.
Bentsen \cite{Bentsen2} (also see \cite{Bentsen}) extended  the result to a class of radon transforms
with fold and blowdown singularities.

However, in higher dimensions ($d\ge 3$) the sharp $L^p$ regularity estimate for $\mathfrak R$ has remained  open for $1<p<2$. 
Set $p_d=2d/(2d-1)$ and
 \[
\alpha(p)=
\begin{cases}
1-\frac 1p, & \ 1\le p<  p_d,
\\
\ \frac 1{2d},  & \ p_d \le p\le  2.
\end{cases}
\]  
 It is natural to conjecture the following.

\begin{conj} \label{sobolev}  
Let $d\ge3$ and $1< p < 2$. Suppose $\gamma$ is a smooth  nondegenerate curve.
Then, $\mathfrak R$ boundedly maps  $L^p$ to $L_\alpha^p$
for $\alpha \le \alpha(p)$.  
\end{conj}

\noindent
Failure of  $L^p$--$L_\alpha^p$ boundedness for $\alpha> \alpha(p)$ can be shown by a slight modification of the examples in \cite{PS06}.  
(See Proposition \ref{nece} below.)
The following is our main result which verifies the conjecture except for some endpoint cases in every dimension $d\ge 3$.

\begin{thm} \label{thm:improving}
Let $d\ge3$ and $1\le p<p_d$. Suppose $\gamma$ is nondegenerate. 
Then, 
\begin{equation}
\label{A1}
\|\mathfrak R f\|_{L^p_{\alpha}(\mathbb R^{d+1})}\le C\|f\|_{p}
\end{equation}
holds if and only if  $\alpha\le 1-1/p$. 
\end{thm}

When $p \in [p_d,2)$, interpolation  with $L^2$--$L_{1/(2d)}^2$ estimate   yields   \eqref{A1}  for  $\alpha <\alpha(p)$ but 
the estimate \eqref{A1} with the endpoint regularity $\alpha=\alpha(p)$, which looks to be a subtle problem,  remains open. 
By a standard scaling argument (\cite{PS06, PS07}) the result in Theorem \ref{thm:improving} can be extended to the curves of finite type. 

A curve $\gamma: I \mapsto \mathbb R^{d}$ is said to be   of finite  type  if  there is an $L=L(s)$ such that  $\sspan \{ \gamma^{(1)}(s), \dots, \gamma^{(L)}(s)\}=\mathbb R^d$ for each $s\in I$,   and  the smallest of such $L(s)$ is called  the type at $s$. The supremum of the type over $s \in I$ is called the maximal type of $\gamma$ (see, e.g., \cite{PS07, HL}).

\begin{cor}\label{type}
Let $d\ge3$, $1\le p<2$,  and $L>d$. Suppose $\gamma$ is a curve of maximal type $L$.
Then,   $\mathfrak Rf$ is bounded  from $L^p(\mathbb R^{d+1})$ to $L_\alpha^p(\mathbb R^{d+1})$  
for $\alpha \le \min( \alpha(p), 1/(Lp))$ if $p\neq (L+1)/L$ when $L\ge 2d-1$, and
if $p \in (1,p_d) \cup (2d/L,2)$ when $d<L <2d-1$.
\end{cor}

For $p\in [2,\infty]$ it  is easy to show the sharp $L^p$ --$L_{1/(Lp)}^p$ estimate, which can be shown by using the $L^2$--$L_{1/(2L)}^2$ estimate and 
interpolation in a similar manner as above. Corollary \ref{type} and Proposition \ref{nece}  completely settle the problem of the optimal Sobolev regularity estimate for $\mathfrak R$ if $L\ge 2d-1$ when $p \neq (L+1)/L$. 
However, some endpoint cases remain left open
not to mention such estimates for the nondegenerate curve.

In this paper, we  make use of the inductive strategy in the recent work of the authors \cite{KLO}, where smoothing properties of 
the (convolution) averaging operator over curves were studied. Exploiting similarity between $\mathfrak R^\ast f$ and the averaging operator, we 
adapt our previous argument. The main new feature of the current paper is use of the decoupling inequality associated with the conical sets generated by curves  (see 
Definition \ref{lambda} and Theorem \ref{decouplingcurve} below).   Compared with our previous work where 
 the averaging operator was decoupled by a class of symbols  adjusted to short subcurves, our new decoupling inequality allows us to dispense with some technicality due to the 
symbols.  The decoupling inequality can also be used to  simplify the argument in \cite{KLO}.

 \smallskip
\noindent{\bf Organization.}
In Section \ref{sec2}, we reduce the proof of Theorem \ref{thm:improving} to obtaining Proposition \ref{newprop}.
We prove a decoupling inequality associated to a nondegenerate curve (Theorem \ref{decouplingcurve}) in Section \ref{sec:dec} which is crucial for the proof of Proposition \ref{newprop}.
The proofs of Proposition \ref{newprop} and Theorem \ref{thm:improving} are given in   Section \ref{sec:it} and Section \ref{order-up}, respectively.
We discuss the sharpness of the smoothing order $\alpha$ in Section \ref{order-up}.

\smallskip

\noindent{\bf Notation.}
For positive constants $A, D$, we denote  $A\lesssim D$ if there exists a (independent)  constant $C$ such that $A\le CD$, where the constant  $C$ may vary from line to line depending on the context.

\section{Estimates with localized frequency}\label{sec2}
In this section, we reduce the proof of Theorem \ref{thm:improving} to showing an inductive  statement (see Proposition \ref{newprop} below).
Afterwards, we obtain some preliminary results which are needed  to prove Proposition \ref{newprop}.

Let us consider the operator 
\[\mathcal R f(x,t)=\chi(t)\int f(x-t\gamma(s),s)\psi(s)ds,\]
which is the dual operator of $\mathfrak R$.
By duality the estimate \eqref{A1} is equivalent to 
\begin{align}\label{dual2}
\| \mathcal Rf\|_{L^p(\mathbb R^{d+1})} \lesssim \|f\|_{L_{-1/p}^p}, \quad 2d<p<\infty.
\end{align}

For the purpose, we closely follow the line of arguments in our previous paper \cite{KLO}.  So, there is  a significant overlap between the current paper and  \cite{KLO}. 
This can be avoided by omitting some shared details. However, we decide to include  them so that  the paper is self-contained and more easily  accessible. 

\subsection{Frequency localized estimate}\label{Sec2.1}
We begin with defining a class of curves in order to prove  \eqref{dual2} in an inductive manner.
For an integer $1 \le L \le d$, by $ \Vol (v_1,\dots,v_L)$ we denote the $L$-dimensional  volume of the parallelepiped generated by 
vectors $v_1,\dots, v_L \in \mathbb R^d$. 

\begin{defn}
Let $B\ge 1$.    We say $\gamma\in \mathfrak V^d(L,B)$ if  $\gamma \in C^{3d+1}(I)$ satisfies
\begin{align}
\label{curveB}
&\max_{s\in I} |\gamma^{(j)}(s)|\,\le\, B,  
   \qquad 0\le j\le 3d+1,
 \\
\label{lindepN} &\min_{s\in I} \Vol \big( \gamma^{(1)}(s),\dots, \gamma^{(L)}(s) \big) \ge B^{-1}. 
\end{align}
\end{defn}

For a smooth function $a(s,t,\xi)$  on $I \times [1,2] \times \mathbb R^d$,
we define 
\[
\mathcal R[a]f(x,t)=(2\pi)^{-d}\iint  e^{i(x-t\gamma(s))\cdot\xi}a(s,t,\xi)  \,\mathcal F_x f (\xi,s)ds d\xi.
\]
Here $\mathcal F_x$ denotes  Fourier transform in $x$. Note that $\mathcal Rf =\mathcal R[a]f$ if $a(s,t,\xi)=\psi(s)\chi(t)$. 
We prove the estimate \eqref{dual2} by induction on $L$ for $\gamma \in \mathfrak V^d(L,B)$
under the localized nondegeneracy assumption:
\begin{align}
\label{sumN}  \sum_{\ell=1}^L & |\langle \gamma^{(\ell)}(s),  \xi\rangle| \ge B^{-1}|\xi|
\end{align}
which holds if  $(s,t,\xi) \in \supp a$ for some $t$. 
When $L<d$,  \eqref{sumN}  can  not be true in general even if $\gamma$ is nondegenerate. However,
an appropriate decomposition in the frequency domain makes it possible that  \eqref{sumN} holds. 
To do this, we consider a class of symbols $a$. 

\begin{defn}\label{def-sym}
Let $\mathbb A_k=\{ \xi \in \mathbb R^d: 2^{k-1} \le |\xi| \le 2^{k+1}\}$ for $k\ge 0$, and $\mathcal I_L=\{(j,\alpha): 0 \le j \le 2L, ~|\alpha| \le d+L+2\}$. 
We say a symbol $a\in C^{d+L+2}(\mathbb R^{d+2})$ is of type $(2^k,L,B)$ if 
$\supp a \subset I \times [2^{-1},2^2]\times \mathbb A_k$,
\begin{align*}
|\partial_t^{j}\partial_\xi^\alpha a(s,t,\xi)| &\le  B |\xi|^{-|\alpha|}, \quad  (j,\alpha) \in \mathcal I_L, 
\end{align*}
 and \eqref{sumN} holds on $\supp_{s,\xi} a$.  Here, as in \cite{KLO}, we denote  $\supp_{s,\xi} a=\cup_{t} \supp a(\cdot,t,\cdot)$. We simply say  \emph{a  statement $S(s,\xi)$, depending on  $s,\xi$,  holds on   $\supp a$ if $S(s,\xi)$ holds for  $s,\xi\in \supp_{a,\xi} a$}.   We also use the same  convention with other variables.
\end{defn}

The estimate \eqref{dual2}  (and hence Theorem \ref{thm:improving}) follows from the next theorem via a standard argument using Fefferman-Stein $\#$-function.
See Section \ref{sec:FS} for details.  

\begin{thm}\label{lclb}
Suppose that $\gamma\in \mathfrak V^d(L,B)$ and $a$ is a symbol of type $(2^k,L,B)$.
Then, for $p> 2L$
\begin{align}\label{smoothing}
\|
\mathcal R[a]f\|_{L^p(\mathbb R^{d+1})}\le
C 2^{-\frac{k}p}\|f\|_p. 
\end{align}
\end{thm}

As we mentioned above, we prove Theorem \ref{lclb} by induction on $L$. 
Theorem \ref{lclb} with $L=1$ is easy to prove. 
Indeed, setting $\tilde{\mathcal R}f=\mathcal F_x(\mathcal R[a] \mathcal F_x^{-1}\!f)$, we note that
\[\tilde{\mathcal R}^*\tilde{\mathcal R}f(\xi,s)=\int \mathcal K(s,s',\xi)f(\xi,s')\,ds',\] 
where 
\[\mathcal K(s,s',\xi)=\int e^{it(\gamma(s)-\gamma(s'))\cdot \xi} \overline a(s,t,\xi) a(s',t,\xi)\,dt.\]
Since \eqref{sumN} holds with $L=1$ on $\supp a$, integration by parts gives
$
|\mathcal K(s,s',\xi)| \le C (1+2^k|s-s'|)^{-2}.
$
By Young's convolution inequality it follows that  $\|\tilde{\mathcal R}^* \tilde{\mathcal R}f \|_2 \lesssim 2^{-k}\|f\|_2$.  
Thus, we get  $\| \mathcal R[a]f\|_2 \lesssim 2^{-k/2}\|f\|_2$  by Plancherel's theorem.
Interpolation with the trivial estimate $\| \mathcal R[a]f\|_\infty \lesssim \|f\|_\infty$ gives \eqref{smoothing} with $L=1$.

Consequently, Theorem \ref{lclb} for $L\ge2$ follows from the next proposition (cf. \cite[Proposition 2.3]{KLO}). 
\begin{prop}\label{newprop}
Let $2 \le N  \le d$. Suppose Theorem \ref{lclb} holds with $L=N-1$.
Then, Theorem \ref{lclb} holds with $L=N$.
\end{prop}

We prove the proposition through the rest of this section, Section \ref{sec:dec} and \ref{sec:it}. Fixing $2\le N\le d$, we  assume that
Theorem \ref{lclb} holds with $L=N-1$.  Additionally,  assuming  that 
 $\gamma \in \mathfrak V^d(N,B)$ and $a$ is of type $(2^k,N,B)$, we prove \eqref{smoothing} for $p>2N$.  For the purpose, composing the symbol $a$, we may further assume that
\begin{align}\label{lowerN}
|\gamma^{(N)}(s)\cdot\xi|\ge (2B)^{-1}|\xi|
\end{align}
holds on $\supp a$.
Otherwise, \eqref{sumN} holds with $L=N-1$, so the  hypothesis  (Theorem \ref{lclb} with $L=N-1$) yields \eqref{smoothing} for $p> 2(N-1)$. 

We prove Proposition \ref{newprop}  in Section \ref{sec:it} using the associated decoupling inequality which is obtained in Section \ref{sec:dec}.
The rest of the section is devoted to proving  two lemmas (Lemma \ref{kernel} and \ref{lem:res}) which play  crucial roles in proving Proposition \ref{newprop}.

\subsection{Symbols adapted to $\gamma$}
We define a class of symbols adapted to the curve $\gamma$.  
From now on,  we assume that $\delta$ satisfies
\begin{align}
\label{delcon}
2^{-k /N}\le \delta \le  (2^2B)^{-N}.  
\end{align}
Let  $\gamma$ satisfy \eqref{lindepN} with $L=N-1$. For $s\in I$, set $\mathrm V_s^{\gamma,\ell} = \sspan \big\{ \gamma^{(j)}(s) : j=1,\dots,\ell\big\}$. Consider  a linear map $\tilde {\mathcal L}_s^\delta :\mathbb R^d \mapsto \mathbb R^d$  given  as follow: 
\begin{equation*}
\begin{aligned}
(\tilde{\mathcal L}_{s}^\delta)^\intercal \gamma^{(j)}(s)
&=\delta^{N-j}\gamma^{(j)}(s), \qquad  && j=1,\dots,N-1,
\\
(\tilde{\mathcal L}_{s}^\delta)^\intercal v &=v, \qquad \qquad   \quad~ && v\in  \big(\mathrm V_s^{\gamma,N-1}\big)^{\perp}.
\end{aligned}
\end{equation*} 
We also consider a linear map $\mathcal L_s^\delta: \mathbb R^{d+1} \mapsto \mathbb R^{d+1}$  given by
\[
\mathcal L_{s}^\delta(\tau, \xi)=
\big(\delta^N \tau - \gamma(s)\cdot   \tilde{\mathcal L}_s^\delta \xi,\, \tilde{\mathcal L}_s^\delta \xi\big), \quad (\tau,\xi)\in \mathbb R \times \mathbb R^d.
\]
Denoting  $G(s)=(1,\gamma(s))$, we set
\begin{align*}
\Lambda_k(\delta,s)=\bigcap_{0\le j\le N-1}\big\{(\tau,\xi) \in \mathbb R \times \mathbb A_k : 
&\,|\langle G^{(j)}(s),(\tau,\xi)\rangle| \le  B2^{k+5} \delta^{N-j}\big\},
\end{align*}
which roughly corresponds to the Fourier support of the operator $\mathcal R[a]f$ with $\supp_s a$ included in an interval centered at $s$ of  length about  $\delta$.  
We define a class of symbols  associated with   $\Lambda_k(\delta,s)$

\begin{defn} 
\label{lambda}
Let $\sz \in (-1,1)$ and $0<\delta \le1$ such that
$I(\sz,\delta):=[\sz-\delta,  \sz+\delta]\subset I$. 
We denote by $\mathfrak A_k(\delta,\sz)=\mathfrak A_k(\delta, \sz,d,N,B,\gamma)$  
the set of smooth functions $\mathfrak a$ on $\mathbb R^{d+3}$ which satisfy the following\,$:$
\begin{align}\label{symsupp}
& \ \ \supp \mathfrak a   \subset I(\sz,\delta)\times [1,2]\times \Lambda_k(\delta,\sz), \\[1ex]
\label{symineq2}
&\big|\partial^{j}_t\partial^{\alpha}_{\tau, \xi}
\mathfrak a\big( s,t,  \mathcal L_{\sz}^\delta(\tau,\xi) \big) \big|
\le B|(\tau,\xi)|^{-|\alpha|}, \quad (j,\alpha ) \in \mathcal I_N. 
\end{align}
\end{defn}

It should be noted that there is  no $s$-differentiation in \eqref{symineq2}.  
Here, $\mathcal I_N$ is given in Definition \ref{def-sym}. We set
\begin{align}\label{T}
\mathcal F(\mathcal T [\mathfrak a] f)(\xi,\tau) =\iint e^{-it'(\tau+ \gamma(s)\cdot \xi)}\mathfrak a(s,t',\tau,\xi)dt' \, \mathcal F_x f (\xi,s)ds.
\end{align}
Clearly, $\mathcal R[a]f=\mathcal T[a]f$ if $\mathfrak a=a(s,t,\xi)$. 
The following is an analogue of \cite[Lemma 2.7]{KLO}.

\begin{lem}\label{kernel}
Let  $\tilde \chi \in \mathrm C_0^\infty((2^{-2},2^2))$ such that $\tilde \chi=1$ on $[3^{-1},3]$.
Let $\mathfrak a$ be a smooth function  which satisfies  \eqref{symsupp} 
and \eqref{symineq2} with $j\le 2$ and $|\alpha| \le d+3$. Then, we have 
\begin{align} 
\label{ker-est} 
\| \mathcal T[\mathfrak a]f\|_{L^p(\mathbb R^{d+1})}  &\le C\delta^{1-\frac 1p} \|f\|_p
\end{align}
for $p\ge2$, and 
\begin{align}
\label{ker-est2}
\|(1-\tilde\chi(t)) \mathcal T[\mathfrak a]f\|_{L^p(\mathbb R^{d+1})} &\le C\delta^{1-\frac 1p-N}2^{-k} \|f\|_p, \quad p>1. 
\end{align}
\end{lem}

\begin{proof}
Note that $\mathcal T[\mathfrak a]f (x,t) = \int K[\mathfrak a](s,t, \cdot)\ast f(\cdot,s)(x) \,ds$
where
\[
K[\mathfrak a](s,t,x)=\frac1{ (2\pi)^{d+1}} \iiint e^{i(t-t')\tau+i(x-t'\gamma(s))\cdot \xi}\mathfrak a(s,t',\tau,\xi) \, d\xi d\tau dt'.
\]
It is easy to show that $|(\mathcal L_{\sz}^\delta)^{-1} \mathcal L_s^\delta(\tau,\xi)| \sim |(\tau,\xi)|$
provided $|s-\sz| \le\delta$ (cf. \cite[Lemma 2.6]{KLO}).
Since \eqref{symineq2} holds with $j=0$ and $|\alpha| \le d+3$, 
it follows that $\supp \mathfrak a(s,t,2^k\mathcal L_s^\delta \cdot) \subset \{(\tau,\xi):|(\tau,\xi)| \lesssim 1\}$
and $|\partial_{\tau,\xi}^\alpha \big( \mathfrak a(s,t,2^k\mathcal L_s^\delta(\tau,\xi))\big)| \lesssim 1$, $|\alpha| \le d+3$.
By changing variables $(\tau,\xi) \rightarrow 2^k\mathcal L_{s}^\delta(\tau,\xi)$ followed by  repeated integration by parts, we have 
\[|K[\mathfrak a](s,t,x)|\lesssim \delta^{\frac {N(N+1)}2}2^{k(d+1)}\int_1^2\big(1+2^k|(\delta^N(t-t'),(\tilde{\mathcal L}_{s}^\delta)^\intercal(x-t\gamma(s)))|\big)^{-d-3}dt'.
\]
This gives  $\|K[\mathfrak a](s,t,\cdot)\|_{L^1_x}\lesssim 1$.  From \eqref{symsupp}, note $\mathcal T[\mathfrak a]f (x,t) = \int_{I(\sz,\delta)} K[\mathfrak a](s,t, \cdot)\ast f(\cdot,s)(x) \,ds$. Thus, we get 
\[\| \mathcal T[\mathfrak a]f\|_{L^\infty(\mathbb R^{d+1})}  \le C\delta \|f\|_\infty.\]
Recall   \eqref{T}.  By translation $\tau \rightarrow \tau-\gamma(s)\cdot \xi$, integration by parts in $t'$, we see
$| \mathcal T[\mathfrak a]f(\xi,\tau)| \lesssim  \int (1+|\tau|)^{-1} |\mathcal F_xf(\xi,s)|\,ds $.
Thus from Plancherel's Theorem and H\"older's inequality, we obtain
\[
\| \mathcal T [\mathfrak a] f\|_2^2
\lesssim  \delta \int_{I(\sz,\delta)}\, \| \mathcal F_x f(\cdot, s) \|_2^2 \,ds  \lesssim \delta \|f\|_2^2.
\]
Therefore, interpolation gives \eqref{ker-est}.  To show \eqref{ker-est2},  we note from the above   estimate for  
$K[\mathfrak a](s,t,x)$ that $\|(1-\tilde\chi(t))K[\mathfrak a](s,t,\cdot)\|_{L^1_x}\lesssim \mathfrak K(t)=: 2^{-k} \delta^{-N} |t-1|^{-1}(1-\tilde\chi(t))$. By  \eqref{symsupp}, using
H\"older's and  Young's convolution inequalities, as before,  we see that $\|(1-\tilde\chi)\mathcal T[\mathfrak a]f\|_p^p$ is bounded above by a constant times 
\begin{align*}
\delta^{p-1}  \int  \mathfrak K^p(t) \int_{I(\sz,\delta)}  \|f(\cdot,s)\|_{L^p_x}^pdsdt
 \lesssim  C\delta^{p-1-pN}2^{-pk} \|f\|_p^p.   
\end{align*}
This gives \eqref{ker-est2}.
\end{proof}

\subsection{Rescaling}
Let $I(\sz,\delta)\subset I$. For  $\gamma\in \mathfrak V^d(N,B)$ we consider a rescaled curve 
\[
\gamma_{\sz}^\delta(s):=
\delta^{-N}(\tilde {\mathcal L}_{\sz}^\delta)^\intercal \big(\gamma(\delta s+\sz)-\gamma(\sz)\big).
\]

\begin{lem}
\label{sc}  Let  $\gamma \in \mathfrak V^d(N,B)$. If  $0<\delta<\delta_\ast$ for a  $\delta_\ast$ small enough, $\gamma_{\sz}^\delta \in \mathfrak V^d(N,3 B)$ and  $\gamma_{\sz}^\delta \in \mathfrak V^d(N-1,B')$ for some $B'$. 
\end{lem}

\begin{proof}
Taylor series expansion of $\gamma^{(j)}(\delta s+\sz)$ at $s=0$ yields
\[
(\gamma_{\sz}^\delta)^{(\ell)}(s)=\sum_{0\le j \le N-1-\ell} \gamma^{(\ell+j)}(\sz)\frac{s^j}{j!}+(\tilde{\mathcal L}_{\sz}^\delta )^{\intercal} \gamma^{(N)}(\sz) \frac{s^{N-\ell}}{(N-\ell)!} +O(B\delta)
\]
for $1 \le \ell \le N-1$ and $(\gamma_{\sz}^\delta)^{(N)}(s)=(\tilde{\mathcal L}_{\sz}^\delta)^{\intercal} \gamma^{(N)}(\sz)+O(B\delta)$. 
Writing $\gamma^{(N)}(\sz)=v_1+v_2 \in \mathrm V_{\sz}^{\gamma,N-1} \oplus (\mathrm V_{\sz}^{\gamma,N-1})^\perp$, we have  $(\tilde{\mathcal L}_{\sz}^\delta)^{\intercal} \gamma^{(N)}(\sz)=(\tilde{\mathcal L}_{\sz}^\delta)^{\intercal} v_1+v_2=v_2+O(B\delta)$. Since $\gamma \in \mathfrak V^d(N,B)$, we see $\gamma_{\sz}^\delta \in \mathfrak V^d(N,3B)$ if $0<\delta<\delta_\ast$ for a sufficiently small $\delta_\ast>0$. In a similar manner, one can also see that $\gamma_{\sz}^\delta \in \mathfrak V^d(N-1,B')$ for some $B'$. 
\end{proof} 

The following lemma, which is  an analogue of \cite[Lemma 2.8]{KLO}, is  important for  our inductive argument.  Let us set
\[
\mathcal  R[\gamma_{\sz}^{\delta},a]f(x,t)=(2\pi)^{-d}\iint  e^{i(x-t\gamma_{\sz}^{\delta}(s))\cdot\xi}a(s,t,\xi)  \mathcal F_x f (\xi,s)\,ds d\xi.
\]

\begin{lem}\label{lem:res}
Let  $\sz \in (-1,1)$, $\mathfrak a \in \mathfrak A_k(\delta,\sz)$, and  $\gamma \in \mathfrak V^d(N,B)$.
Suppose   
\begin{align}\label{lowscale}
\textstyle \sum_{j=1}^{N-1}\, \delta^{j} |\langle \gamma^{(j)}(s), \xi \rangle |\ge
B^{-1}2^k \delta^N 
\end{align}
for $(s,\xi)\in I(\sz,\delta)\times \supp_\xi \fa$.  Then, there exist constants  $C$, $\tilde B$, $\delta_{*}=\delta_{*}(B,N,d)$, and   
$\tilde f$ and  a symbol $\tilde a$ such that 
\begin{align}\label{resc}
\| \tilde \chi(t) \mathcal T[\mathfrak a]f \|_p=\delta^{1-\frac 1p} \| \mathcal R[\gamma_{s_0}^\delta,\tilde a]\tilde f\|_p
\end{align}
for  $0<\delta<\delta_{*}$,  $\|\tilde f\|_p=\|f\|_p$, 
$|\partial_t^j\partial_\xi^\alpha \tilde a(s,t,\xi)| \le \tilde B|\xi|^{-|\alpha|}$ for $(j,\alpha)\in \mathcal I_{N-1}$, 
and 
\Be
\label{asupp}
\supp \tilde a_\xi \subset I\times [2^{-2},2^2]\times \{\xi \in \mathbb R^d: C^{-1}\delta^N2^k\le |\xi| \le C \delta^N2^k\}.
\Ee
\end{lem}

\begin{proof}
Let  $\mathfrak a_\delta(s,t,\tau,\xi)=\mathfrak a(\delta s+\sz,t,\tau,\xi)$.
By Fourier inversion and \eqref{T},  changing  variables $s\rightarrow \delta s+\sz$, $(\tau,\xi)\rightarrow (\tau-\gamma(\sz)\cdot \xi, \xi)$ gives
\Be\label{ta}
\mathcal T[\mathfrak a]f(x,t)=(2\pi)^{-d}\delta \iint e^{i\langle x-t\gamma(\sz),\xi \rangle} b(s,t,\xi) \mathcal F_xf(\xi, \delta s+\sz)\,ds d\xi,
\Ee
where
\[
b(s,t,\xi)= \frac{1}{2\pi}\iint e^{it\tau} e^{-it' ( \tau+ \langle \gamma(\delta s+\sz)-\gamma(\sz), \xi \rangle)} \mathfrak a_\delta\big(s, t', \tau-\gamma(\sz)\cdot \xi, \xi\big)\,dt' d\tau.
\]
We observe that
\[\tilde \chi(t) b(s,t, \delta^{-N}\tilde{\mathcal L}_{\sz}^\delta \xi) = e^{-it\gamma_{\sz}^\delta (s)\cdot \xi}\, \tilde a(s,t,\xi),\]
where
\begin{align}\label{atilde}
\tilde{ a}(s,t,\xi)&=\frac 1{2\pi}\iint e^{-it'(\tau+\gamma_{\sz}^\delta(s)\cdot \xi)}\tilde{\chi}(t)\mathfrak a_\delta(s,t'+t,\delta^{-N}\mathcal L_{\sz}^\delta(\tau,\xi))\,dt'd\tau.
\end{align}
It is clear that   \eqref{asupp} holds for some $C\ge1$.
Since $\mathfrak a \in \mathfrak A_k(\delta,\sz)$, it is not difficult to see 
$
|\partial_t^j \partial_\xi^\alpha \tilde a(s,t,\xi)| \le \tilde B |\xi|^{-|\alpha|}$ for  $(j,\alpha) \in \mathcal I_{N-1}$
(see (2.25) in \cite{KLO}).

Set $\mathcal C_p=\mathcal C_p(\delta):= \delta^{1/p} |\det \delta^{-N} \tilde{\mathcal L}_{\sz}^\delta|^{1-1/p}$. Let $\tilde f$ be given by  $\mathcal F_x \tilde f (\xi,s)= \mathcal C_p \mathcal F_x f( \delta^{-N} \tilde{\mathcal L}_{\sz}^\delta \xi, \delta s+\sz)$, thus $\|\tilde f\|_p=\|f\|_p$.  Recalling \eqref{ta} and changing variables  $\xi \rightarrow \delta^{-N} \tilde{\mathcal L}_{\sz}^\delta \xi$, we now have 
\[
\tilde \chi(t) \mathcal T[\mathfrak a]f(x,t) =  \frac{\mathcal C_{p'}}{(2\pi)^{d}}
\iint e^{i \langle x-t\gamma(\sz), \delta^{-N}\tilde{\mathcal L}_{\sz}^\delta \xi \rangle} e^{-it\gamma_{\sz}^\delta (s)\cdot \xi} \,\tilde a(s,t,\xi) \mathcal F_x \tilde f(\xi,s)\,ds d\xi.
\]
This gives $\tilde \chi(t)\mathcal T [\fa]f(x,t)
=\mathcal C_{p'} \mathcal R[\gamma_{\sz}^\delta,\tilde a]
\tilde f \big( y,t\big)$ 
where $y=\delta^{-N} ( \tilde{\mathcal L}_{\sz}^\delta)^{\intercal} (x-t\gamma(\sz))$. 
Therefore, changing variable $x \rightarrow \delta^{N} (\tilde{\mathcal L}_{\sz}^\delta)^{-\intercal}x+t\gamma(\sz)$, we obtain \eqref{resc}.
\end{proof}

Combining Lemma \ref{lem:res} and the hypothesis (Theorem \ref{lclb} with $L=N-1$), we obtain  the following.

\begin{cor}\label{lessdegen}
Suppose that Theorem \ref{lclb} holds with $L=N-1$,  and  $\mathfrak a$,  $\gamma$, and $\delta_\ast$ are  the same as in Lemma \ref{lem:res}.  
Then, 
if $p> 2(N-1)$,  for $0<\delta<\delta_*$ we have
\[
\big\|   \mathcal T [\fa] f \big\|_p
\lesssim   2^{-\frac kp}\delta^{1-\frac {N+1}p} \|f\|_p.
\]
\end{cor}

\begin{proof}
By \eqref{resc} and dyadic decomposition (of $\tilde a$ in the Fourier side), we have
\begin{align}\label{res:TA}
\big\|  \tilde\chi \mathcal T [\fa] f \big\|_p
\le C  \delta^{1-\frac 1p} \sum_{0\le \ell \le C }
\big\|\mathcal  R[\gamma_{\sz}^{\delta},a_\ell]
f_\ell \big\|_p,  
\end{align}
for some constant $C$  where $\|f_\ell\|_p=\|f\|_p$, and $a_\ell$ are symbols of type $(2^j,N-1, \tilde B)$ with $C^{-1}2^k\delta^N\le 2^j \le C2^k\delta^N$. Once we have this, the proof is straightforward. 
By Lemma \ref{sc}, $\gamma_{\sz}^\delta \in \mathfrak V^d(N-1, B')$ for some $B'>0$. 
Since $\|f_l\|_p =\|f\|_p$, applying Theorem \ref{lclb} with $L=N-1$, we have
\[
\textstyle \big\|  \tilde\chi \mathcal T [\fa] f \big\|_p \le C \sum_{l}\delta^{1-\frac 1p} (2^k\delta^N)^{-\frac 1p} \|f_l\|_p  \lesssim  2^{-\frac kp}\delta^{1-\frac {N+1}p} \|f\|_p
\]
for $ p> 2(N-1)$. Recalling \eqref{delcon}, we combine this and \eqref{ker-est2} to get the desired bound.

It remains to show \eqref{res:TA}. In fact, after applying Lemma \ref{lem:res}  we only need to adjust the support of the consequent symbol $\tilde a$ via by moderate decomposition and scaling.
We omit details. (See the proof of \cite[Lemma 2.8]{KLO}.) 
\end{proof}

\section{Decoupling inequality for curve}\label{sec:dec}

In this section, we prove the decoupling inequality, which is to be used to decompose the operator $\mathcal T[\mathfrak a]f$.  
In our earlier work \cite{KLO}, the averaging operator was decoupled by making use of decomposition based on a class of symbols which are adjusted to short subcurves. 
The same approach also works to prove  Proposition \ref{newprop}. However, instead of following the previous strategy, we directly obtain a decoupling inequality associated with the conic sets 
\[\Lambda_k(\delta,s_l), \quad 1\le l\le L,\] while 
$\{s_1,\dots,s_L\}\subset I$ is a collection of $\delta$-separated points contained in $I$. More precisely,  we have the following.

\begin{thm}\label{decouplingcurve}
Let $0<\delta\le1$ and $S:=\{s_1,\dots,s_L\}\subset I$  be a collection of $\delta$-separated points.
Then, if $2\le p\le N(N+1)$, for  any $\epsilon>0$  there is a constant $C_\epsilon=C_\epsilon(B)$, independent of $S$,  such that
\begin{align}
\label{eq-decouple}
\big\| \sum_{1 \le l\le L} f_l \big\|_{L^p(\mathbb R^{d+1})}\le C_\epsilon   \delta^{-\epsilon}
\Big(\sum_{1 \le l\le L} \|f_l\|_{L^p(\mathbb R^{d+1})}^2
\Big)^{1/2}
\end{align} 
holds whenever $\supp \widehat {f}_{l}\subset \Lambda_k(\delta,s_l)$.
\end{thm}

H\"older's inequality gives $\big\| \sum_{1 \le l\le L} f_l \big\|_{p}\le C_\epsilon   \delta^{-\epsilon} \delta^{1/p-1/2} 
(\sum_{1 \le l\le L} \|f_l\|_{p}^p)^{1/p}$. Interpolation with the trivial $L^\infty-\ell^\infty L^\infty$ estimate yields the inequality 
\begin{align}\label{decoupta}
\big\| \sum_{1 \le l\le L} f_l \big\|_{L^p(\mathbb R^{d+1})}
	\le C_{\epsilon}
	\delta^{-1+\frac {N+1} {p}+\epsilon} \Big(\sum_{1 \le l\le L} 
	\big\| f_l \big\|_{L^p(\mathbb R^{d+1})}^p \Big)^{\frac 1p}
\end{align}
for $p>2N$ whenever $\supp \widehat {f}_{l}\subset \Lambda_k(\delta,s_l)$. 

\subsection{Decoupling inequality for curve}
Fixing $N\ge 2$, we now consider the slabs given by  an anisotropic neighborhood of the moment curve 
\[
\gamma_{\circ}(s):=\big(s,\ {s^2}/{2!},\dots, {s^{N+1}}/{(N+1)!}\big).
\]

\begin{defn}  Let $0<\delta\le1$ and $B\ge1$.  For $s \in I$, let $\mathbf S(s,\delta, B)$ denote the set of $(\tau,\xi) \in \mathbb R \times \mathbb R^N$ such that
\begin{align*}
B^{-1} \le | \langle \gamma_{\circ \phantom{1}}^{(N+1)}(s), (\tau,\xi) \rangle|  \le B; \quad  
| \langle \gamma_{\circ \phantom{1}}^{(j)}(s),(\tau,\xi) \rangle|  \le   \delta^{N+1-j}, \ \ j=1,\dots,N.
\end{align*}
\end{defn}

We now recall the decoupling inequality for  such slabs as above  which was shown in \cite{BGHS} (see also \cite[Corollary 2.15]{KLO}).

\begin{thm}\label{rcdecoupling}
Let $0<\delta\le1$ and $\{s_1, \dots, s_L\}\subset I$ be a collection of $\delta$-separated points contained in $I$.
Denote $\mathbf S_l=\mathbf S(s_l,\delta, B)$.
Then, if   $2 \le p\le N(N+1)$,  for any $\epsilon>0$ there is a constant $C_\epsilon=C_\epsilon(B)$ such that
\begin{align*}
\big\| \sum_{1 \le l\le L} F_l \big\|_{L^p(\mathbb R^{N+1})}\le C_\epsilon  \delta^{-\epsilon}
\Big(\sum_{1 \le l\le L} \|F_l\|_{L^p(\mathbb R^{N+1})}^2
\Big)^{1/2}
\end{align*} 
holds whenever $\supp \widehat {F}_{l}\subset \mathbf S_l$. 
\end{thm}

To show Theorem \ref{decouplingcurve}, we apply the decoupling inequality after projecting the sets $\Lambda_0(\delta,s_l)$ to the subspace $\mathrm V_\mu$ which is spanned by
$\{G^{(0)}(\mu) ,\dots, G^{(N)}(\mu)\}$. To do so, 
for $\mu\in I$ we consider a coordinate system  $\mathbf y_\mu=\mathbf y_\mu(\tau,\xi)$ given by 
\begin{align}\label{ymu}
\mathbf y_\mu=(y_\mu^0,\cdots, y_\mu^N)=(\langle G^{(0)}(\mu ),(\tau,\xi)\rangle, \dots, \langle G^{(N)}(\mu ),(\tau,\xi)\rangle).
\end{align}
Recall that $\gamma\in \mathfrak V^d(N,B)$, so $\Vol(\langle G^{(0)}(\mu ), \dots, G^{(N)}(\mu ))\ge 1/B$.  
Let $\delta, \delta'$ be positive numbers satisfying
\begin{equation}\label{rel}
0< \delta<\delta' \le \delta^{N/(N+1)} \le1.
\end{equation}
Then it is easy to see that
\begin{align}\label{ordering}
(\delta')^{\ell+1} \le \delta^{\ell}, \quad \ell=1,\dots,N.
\end{align}

 The following lemma shows that the projections of the sets $\Lambda_0(\delta,s_l)$ form
a reverse $\delta/\delta'$-adapted cover after a proper linear change of variables (cf. \cite[Lemma 3.3]{KLO}) if $s_l$ are contained in an interval of length $\delta'$. 
Let $\mathrm D_\delta$ denote the $(N+1)\times (N+1)$ diagonal matrix given by 
\[ \mathrm D_\delta=(\delta^{-N}e_1,\delta^{-N+1}e_2,\dots, \delta^{0}e_{N+1}).\]

\begin{lem}\label{scaled}
Let  $\delta, \delta'$ be positive numbers satisfying \eqref{rel} and $s' \in [\mu-\delta',\mu+\delta']$.
Suppose $(\tau,\xi) \in \Lambda_0(\delta,s')$.   Then we have 
\begin{align}
\qquad (4B)^{-1} \le |\langle \mathrm D_{\delta'} \mathbf y_\mu,\gamma_{\vphantom{1}\circ}^{(N+1)} &
 \rangle | \le 4B,  
\label{claim2}
\\
\big|\big\langle \mathrm D_{\delta'} \mathbf y_\mu,\gamma_{\vphantom{1}\circ}^{(j)}\big( \frac{s'-\mu}{\delta'}\big) \big\rangle\big|
	&\lesssim B\big({\delta}/{\delta'}\big)^{N+1-j}, \quad 1 \le j  \le N.  \label{claim1}
\end{align}
\end{lem}

\begin{proof}
Note that \eqref{claim2} is clear from \eqref{lowerN}. To prove \eqref{claim1}, 
we first note that $\langle \mathbf y_\mu,\gamma_{\vphantom{1}\circ}^{(j)}(s) \rangle=(\delta')^{N+1-j} \langle \mathrm D_{\delta'} \mathbf y_\mu, \gamma_{\vphantom{1}\circ}^{(j)}(s/\delta') \rangle$. Thus, it is sufficient  to show that 
\[|\langle \mathbf y_\mu, \gamma_{\vphantom{1}\circ}^{(j)}(s'-\mu) \rangle| \lesssim  B{\delta}^{N+1-j}\] for $1 \le j  \le N$. 
Recalling \eqref{ymu},    we observe
\[ \langle \mathbf y_\mu, \gamma_{\vphantom{1}\circ}^{(j)}(s'-\mu) \rangle 
=
\Big\langle\sum_{\ell=j-1}^N G^{(\ell)}(\mu)\frac{(s'-\mu)^{\ell-j+1}}{(\ell-j+1)!},(\tau,\xi)\Big\rangle. \]
 Taylor's theorem gives \[ \Big |G^{(j-1)}(s')-\sum_{\ell=j-1}^N G^{(\ell)}(\mu)\frac{(s'-\mu)^{\ell-j+1}}{(\ell-j+1)!}\Big|\le B |s'-\mu|^{N-j+2}\]
for $j=1,\dots,N$.
Since $|s'-\mu|\le \delta'$ and $(\tau,\xi) \in \Lambda_0(\delta,s')$,  \eqref{claim1} follows by \eqref{ordering}. 
\end{proof}

By Lemma \ref{scaled} and Theorem \ref{rcdecoupling}, we can show that \eqref{eq-decouple} holds if a $\delta$-separated set $\{s_1, \dots, s_L\}$ are contained in an interval of length $\lesssim \delta^{N/(N+1)}$.  More precisely, we have the following.

\begin{lem}\label{s-dec}
Let $0<\delta\le1$ and $\delta\le \delta'\le \delta^{N/(N+1)}$. Let $\{s_1 ,\dots,s_L \}\subset [\mu-\delta' ,\mu+\delta' ] $  be a collection of $\delta$-separated points.
Then, if $2\le p\le N(N+1)$, for  any $\epsilon>0$  there is a constant $C_\epsilon=C_\epsilon(B)$  such that
\eqref{eq-decouple} holds whenever $\supp \widehat {f}_{l}\subset \Lambda_k(\delta,s_l )$.
\end{lem}

\begin{proof}
Set $\mathrm V_\mu=\sspan \{ \gamma'(\mu),\dots, \gamma^{(N)}(\mu)\}$ and let  $\{v_{N+1},\dots,v_d\}$ be an orthonormal basis of  $\mathrm V_\mu^\perp$ .
Recalling that \eqref{lindepN} holds with $L=N$, we write
$\xi=\overline \xi+\sum_{j=N+1}^d y_j(\xi)v_j$ for $\overline \xi \in \mathrm V_\mu$.
Changing of variables
\[ (\tau,\xi) \rightarrow \mathrm Y_\mu(\tau,\xi):= (\mathbf y_\mu(\tau,\xi), y_{N+1}(\xi),\dots,y_d(\xi))\] (see \eqref{ymu}),
we may work with  the coordinate system given  by $\{ \mathbf y_\mu,y_{N+1},\dots,y_d\}$
instead of $(\tau,\xi)$. We consider  the linear map  
\[\mathrm Y_\mu^{\delta'} (\tau,\xi)=(\mathrm D_{\delta'} \mathbf y_\mu(\tau,\xi), y_{N+1}(\xi),\dots,y_d(\xi)).\]
Since  $\{s_1 ,\dots,s_L \}\subset [\mu-\delta',\mu+\delta'] $ and $\delta'\le  \delta^{N/(N+1)}$, by Lemma \ref{scaled}  it follows that 
\Be
\label{supp-supp}
\mathrm Y_\mu^{\delta'} ( \Lambda_0(\delta,s_l ) )
\subset  \mathbf S_l:=\mathbf S\Big(\frac{s_l -\mu}{\delta'}, C\frac{\delta}{\delta'}, 4B\Big) \times \mathbb R^{d-N}
\Ee
for some $C>0$ depending only on $B$.  Applying Theorem \ref{rcdecoupling} with $\delta$ replaced by $C\delta/\delta'$ and slabs $\mathbf S_l, 1\le l\le L$,
and then using  a trivial extension via Minkowski's inequality, we have 
\[
\textstyle 
\big\| \sum_{1\le l\le L}f_l \big\|_p
	\le C_{\epsilon}
	\delta^{-\epsilon} \big(\sum_{1\le l\le L}
	\big\| f_l \big\|_p^2 \big)^{1/2}
\]
for $2\le p \le N(N+1)$ whenever $\widehat f_l\subset  \mathbf S_l$. Since the decoupling inequality is invariant under affine changes of variables, by undoing the change of variables $(\tau,\xi) \rightarrow \mathrm Y_\mu(\tau,\xi)$ and rescaling  $(\tau,\xi) \rightarrow  2^{-k}(\tau,\xi)$, we obtain  \eqref{eq-decouple} whenever $\supp \widehat {f}_{l}\subset \Lambda_k(\delta,s_l )$.
 \end{proof}

\subsection{Proof of Theorem \ref{decouplingcurve}}
We now prove Theorem \ref{decouplingcurve}. 
Let $2 \le p\le N(N+1)$.  For the purpose, for some $\alpha>0$ we assume that  
\begin{align}
\label{induct}
\tag*{$\mathfrak D(\alpha)$}
\textstyle \big\| \sum_{1 \le l\le L} f_l \big\|_{L^p(\mathbb R^{d+1})}\le C   \delta^{-\alpha}
\big(\sum_{1 \le l\le L} \|f_l\|_{L^p(\mathbb R^{d+1})}^2
\big)^{1/2}
\end{align} 
holds for  $0<\delta\le \delta_0:= (2^2B)^{-N-1}$ with a constant $C$, independent of $S$,   whenever  $\supp \widehat {f}_{l}\subset \Lambda_k(\delta,s_l)$, $1 \le l\le L$.  
Of course, \ref{induct} holds true if $\alpha \ge 1/2$ by Minkowski's and  H\"older's inequalities.  We set \[
\delta'=\delta^{N/(N+1)}.
\]

Let us denote $I _\nu$, $1\le \nu\le M$,  be disjoint intervals of length $\rho\in (2^{-3}\delta', 2^{-2}\delta']$ which partition $I$. 
Let $s _\nu'$ be a point contained in $I _\nu$ such that  $s_1', \dots, s_M'$ are separated at least  by $ 2^{-4}\delta'$. We now claim  that 
\Be\label{inclusion} \Lambda_k(\delta,s_l )\subset \Lambda_k(\delta',s' _\nu)\Ee
if $s_l\in I _\nu$.  Indeed, by scaling it is sufficient to show $\Lambda_0(\delta,s_l )\subset \Lambda_0(\delta',s' _\nu).$  Let  $(\tau,\xi)\in \Lambda_0(\delta,s_l )$. Then, it follows that   $|\langle G^{(\ell)}(s_l),(\tau,\xi)\rangle|\le 2^5B \delta^{1/(N+1)}(\delta')^{N-\ell}$.  
By Taylor's theorem we have
\[\langle G^{(j)}(s' _\nu),(\tau,\xi)\rangle=\sum_{\ell=j}^{N-1}\langle G^{(\ell)}(s_l),(\tau,\xi)\rangle\frac{(s'_\nu-s_l)^{\ell-j}}{(\ell-j)!}+\mathcal E,\]
where $|\mathcal E|\le 2B|s'_\nu-s_l|^{N-j}$. Therefore,  we see that $(\tau,\xi)\in \Lambda_0(\delta',s'_\nu)$. 

Let $\supp \widehat {f}_{l}\subset \Lambda_k(\delta,s_l)$, $1 \le l\le L$. 
We write $\sum_{1\le l\le L}f_l=\sum_{1\le \nu\le M}\sum_{s_l\in I_\nu} f_l$. By \eqref{inclusion}  the Fourier support  of $\sum_{s_l\in I_\nu} f_l$ is included in  $\Lambda_k(\delta',s' _\nu)$.
Since $s'_\nu$ are separated by $2^{-4}\delta'$, \ref{induct} implies 
\[
\big\| \sum_{1 \le l\le L} f_l \big\|_{L^p(\mathbb R^{d+1})}\le 
C  \delta^{-\frac{N\alpha}{N+1}}
\Big(\sum_{1 \le \nu \le M} \|\sum_{s_l\in I_\nu} f_l\|_{L^p(\mathbb R^{d+1})}^2
\Big)^{1/2}
\]
for a constant $C$. Since the length of interval $I_\nu$ is less than $\delta^{N/(N+1)},$ by Lemma \ref{s-dec} we have 
$ \|\sum_{s_l\in I_\nu} f_l\|_{p}\le C_\epsilon  \delta^{-\epsilon}  (\sum_{s_l\in I_\nu}\| f_l\|_{p}^2
)^{1/2}$.  Therefore, combining this and the above inequality, we obtain 
\[
\textstyle
\big\| \sum_{1 \le l\le L} f_l \big\|_{L^p(\mathbb R^{d+1})}\le 
C_\epsilon  \delta^{-\frac{N\alpha}{N+1}-\epsilon}
\big(\sum_{1 \le l\le L} \|f_l\|_{L^p(\mathbb R^{d+1})}^2
\big)^{1/2}
\]
for a constant $C_\epsilon$. This establishes the implication $\mathfrak D(\alpha)\to  \mathfrak D(\epsilon + {N\alpha}/(N+1))$. 
Iteration of this implication suppresses $\alpha$ arbitrarily small.  \qed

\section{Proof of Proposition \ref{newprop}}\label{sec:it}
In this section, we prove Proposition \ref{newprop} by making use of the decoupling inequality \eqref{decoupta}.
As mentioned in Section \ref{Sec2.1} (below Proposition \ref{newprop}), in order to prove Proposition \ref{newprop},
it suffices to show Theorem \ref{lclb} with $L=N$. We first reduce the matter to obtaining estimates for $\mathcal T[\mathfrak a_0]$ with a suitable $\mathfrak a_0$.

\subsection{Reduction}\label{sec:reduce} We begin  by recalling $\gamma \in \mathfrak V^d(N,B)$ and $a$ is of type $(2^k,N,B)$. 
Let $\delta_*$ be the small number given in Lemma \ref{lem:res} and set
\Be\label{delta-min}
\delc=\min\{ \delta_*, (2^2B)^{-N}\}. 
\Ee
Let $\beta_0\in \mathrm C_0^\infty([-1,1])$ such that $\beta_0= 1$ on $[-1/2,1/2]$. 
We set
\[
a_N(s,t,\xi)=a(s,t,\xi)\prod_{1\le j\le N-1} \beta_0 \Big(100dB2^{-k}\delc^{-N} \langle \gamma^{(j)}(s),\xi\rangle \Big).
\]
Clearly, \eqref{sumN} holds  on  $\supp (a-a_N)$ with $L=N-1$ and  $B$ replaced by $(100dB)^{-1}\delc^N$. Since 
$a$ is of type $(2^k,N,B)$,  it is easy to see 
$(a-a_N)$ is a symbol of type $(2^k,N-1,B')$ for some $B'$. 
Thus,  the  hypothesis  (Theorem \ref{lclb} with $L=N-1$ and $B=B'$) gives the estimate 
\[
\|\mathcal R[a-a_N]f\|_p\lesssim 2^{-\frac kp}\|f\|_p
\]
for $p> 2(N-1)$. 
So, we need only to  consider $\mathcal R[a_N]$ instead of $\mathcal R[a]$.
Furthermore, by a moderate decomposition of $a_N$  we assume  
\[\supp_s a_N \subset [\sz-\delc, \sz+\delc]\] for some $\sz \in (-1,1)$. 
We may assume that $\sz=\delc\nu$ for $\nu \in \mathbb Z$.

It is not difficult to see  that the contribution of the frequency part  $\{(\tau,\xi): |\tau+\gamma(s)\cdot \xi| \gtrsim  2^{k+1} \delc^N, \, \forall s \in I\}$ is not  significant. 
To see this, let us set
\[
{\mathfrak a}_0(s,t,\tau,\xi)= a_N(s,t,\xi) 	\beta_0 \big( 
	\delc^{-2N} 2^{-2k}|\tau+\gamma(s)\cdot \xi |^2\big)
\]
and $\mathfrak a_1=\mathfrak a_0-a_N$. Recalling \eqref{T}, by  Fourier inversion we have
\[
\mathcal R[a_N]f = \mathcal T[\mathfrak a_0]f+ \mathcal T[\mathfrak a_1]f.
\]
The operator $ \mathcal T[\mathfrak a_1]$  is easy to handle.  Let us set 
$\mathfrak a=-i2^{k}\delc^{N}(\tau+\langle \gamma(s), \xi \rangle )^{-1} \partial_{t}\mathfrak a_1$. 
Then, by integration by parts in $t'$ and \eqref{T} we see $\mathcal T [\mathfrak a_1]=(2^{k}\delc^{N})^{-1}\mathcal T [\mathfrak a]$. Note that  $|\tau+\gamma(s)\cdot\xi| \gtrsim  2^{k} \delc^N$ on $\supp \mathfrak a_1$ and so on $\supp \mathfrak a$.  It is clear that 
%
$\mathfrak a$ satisfies  \eqref{symsupp} and \eqref{symineq2} with $\delta=\delc$ and $B=C_1\delc^{-C}$ for some large $C,C_1$. Thus, Lemma \ref{kernel} gives
$
\|\mathcal T [ \mathfrak a_1]f\|_p\lesssim 2^{-k}\|f\|_p
$
for $ p\ge 2. $

Therefore,  the proof of Theorem \ref{lclb} with $L=N$ is now reduced to showing that
\begin{align}\label{overline-a0}
\big\| \mathcal T[\mathfrak a_0]f \big\|_{p} 
\le C 2^{-\frac{k} p}\|f\|_{p}, \quad p>2N. 
\end{align}

\subsection{Decomposition} For $n\ge 0$, let us set $ \delta_n=2^n2^{-k/N}$ and 
\Be\label{Jnmu}
\mathfrak J_n= \delta_n\mathbb Z\cap I.
\Ee
We  consider 
\[
\mathfrak G_N(s,\tau,\xi)=\sum_{0\le j\le N-1}(2^{-k}|\langle G^{(j)}(s),(\tau,\xi)\rangle|)^{\frac{2N!}{N-j}},
\]
by which we can decompose  $\mathfrak a_0$ into the symbols contained in  $\mathfrak A_k(\delta_n,s)$ for $s\in \mathfrak J_n$.

Set $\beta_\ast=\beta_0-\beta_0(2^{2N!}\cdot)$. Note that $\beta_0+\sum_{n\ge 1}  \beta_\ast(2^{-2N! n} \cdot)=1$. 
Let $\zeta \in \mathrm C_0^{\infty}([-1,1])$ such that $\sum_{\nu \in \mathbb Z} \zeta(\cdot-\nu)=1$.
We set
\begin{align*} 
\mathfrak a^{n}_\nu
	=\mathfrak a_0
	\times 
	\begin{cases}
	\beta_0 ( 	\delta_0^{-2N!} \, \mathfrak G_N  ) \, \zeta(\delta_0^{-1}s-\nu),  \quad\, \nu \in \mathfrak J_0, \quad   n=0, 
	\\[2pt]
		\beta_\ast ( \delta_n^{-2N!}
	\, \mathfrak G_N) \,	\zeta(\delta_n^{-1}s-\nu), \quad\,   \nu \in \mathfrak J_n,\quad  n\ge 1.
	\end{cases}
\end{align*}
 Then, it follows that 
\begin{align}\label{split}
\mathfrak a_0(s,t,\tau,\xi)= \sum_{n\ge 0}\sum_{\nu \in \mathfrak J_n} \mathfrak a_\nu^{n}(s,t,\tau,\xi).
\end{align}
Since $\delc$ is the fixed constant, it is clear that $C^{-1} \mathfrak a_0 \in \mathfrak A_k(\delc,\sz)$ for a large constant  $C>0$. So, $\supp \mathfrak a_0 \subset\Lambda_k(\delc,\sz)$ and
$\mathfrak G_N\lesssim 1$ for $(\tau,\xi)\in\supp \mathfrak a_\nu^n$. Obviously, we may assume  $\delta_n\lesssim 1$ since  $\mathfrak a_\nu^n=0$ otherwise.

The following tells that $\mathfrak a_{\nu}^{n}$ is contained in a proper symbol class.

\begin{lem}[cf. {\cite[Lemma 3.2]{KLO}}]\label{symcheck}
For $n\ge0$, there exists a constant $C$ such that $C^{-1}\mathfrak a_{\nu}^{n}\in \mathfrak A_k(\delta_n,\delta_n\nu)$.
\end{lem}
\begin{proof}
The condition \eqref{symsupp} trivially holds for $\mathfrak a=\mathfrak a_\nu^n$. 
So, we only need to show \eqref{symineq2} for $\delta=\delta_n$ and $s=\delta_n\nu$.

It is not difficult to see that $\mathfrak a_0$ satisfies \eqref{symineq2} (see \cite[(3.35)]{KLO}).
So it suffices to show \eqref{symineq2} for $\beta_N(\delta_n^{-2N!}\mathfrak G_N(s,\tau,\xi))$. 
By Leibniz's rule, it is enough to prove that
\begin{align}\label{deriv1}
|\nabla_{\tau,\xi}\delta_n^{-(N-j)} 2^{-k}  \big\langle G^{(j)}(s), {\mathcal L}_{\delta_n\nu}^{\delta_n}(\tau, \xi) \big\rangle|
\lesssim\! 2^{-k},
\end{align}
for $j=0,\dots,d-1$.
Note that if $|\delta_n-s|\le \delta_n$, then
\begin{align}\label{simsim}
|({\mathcal L}_{s}^{\delta_n})^{-1}{\mathcal L}_{\delta_n \nu}^{\delta_n}(\tau,\xi)|\sim |(\tau,\xi)|
\end{align}
(see \cite[Lemma 2.6]{KLO}).
Recall that $\nabla_{\tau,\xi}\langle G^{(j)}(s), {\mathcal L}_{\delta_n\nu}^{\delta_n}(\tau, \xi)\rangle
= ({\mathcal L}_{\delta_n\nu}^{\delta_n})^\intercal G^{(j)}(s)$. Thus, by 
 \eqref{simsim} we get \eqref{deriv1}.
\end{proof}

\subsection{Proof of Proposition \ref{newprop}}
By the reduction in Section \ref{sec:reduce}, it suffices to prove \eqref{overline-a0}.
Recalling \eqref{split} and applying  the Minkowski inequality, we have
\[
\| \mathcal T[\mathfrak a_0] f\|_p
\le \sum_{2^{-k/N} \le \delta_n\lesssim 1}  \big\| \sum_{\nu \in \mathfrak J_n} \mathcal T[\mathfrak a_{\nu}^{n}]f \big\|_p.
\]
Using  Lemma \ref{symcheck}, one can easily see that $\supp \mathfrak a_{\nu}^{n}\subset \Lambda_k(\delta_n, \delta_n\nu ) $.  
Thus, we may use the decoupling inequality  \eqref{decoupta}. 
Combining this  and the above inequalities gives  
\[
\textstyle \| \mathcal T[\mathfrak a_0] f\|_p
	\le
	C_{\epsilon}  \sum_{2^{-k/N} \le \delta_n\lesssim 1} {\delta_n}^{-1+\frac {N+1} {p}+\epsilon}
	\big(  \sum_{\nu \in \mathfrak J_n} \big\| \mathcal T[\mathfrak a_{\nu}^{n}]f \|_p^p \big)^{1/p}
\]
for $2N<p<\infty$. Hence, for the estimate \eqref{overline-a0}  it suffice to show that 
\begin{equation}
\label{tannu}
\|\mathcal T [\mathfrak a_{\nu}^{n}]f \|_p\lesssim \delta_n^{1-\frac{N+1}p}2^{-\frac kp}\|f \|_p, \quad p>2N. 
\end{equation}
Indeed, let $f_\nu(x,s)=\tilde \zeta(\delta_n^{-1}s-\nu) f(x,s)$ where
$\tilde \zeta\in \mathrm C_0^\infty([-2,2])$ such that $\tilde \zeta=1$ on $\supp \zeta$.  From \eqref{T} we see 
$\mathcal T [\mathfrak a_{\nu}^{n}]f=\mathcal T [\mathfrak a_{\nu}^{n}]f_\nu$. Combining this and \eqref{tannu}, we have 
\[
\textstyle \big(  \sum_{\nu \in \mathfrak J_n} \big\| \mathcal T[\mathfrak a_{\nu}^{n}]f \|_p^p \big)^{1/p}\lesssim\delta_n^{1-\frac{N+1}p}2^{-\frac kp}\big(  \sum_{\nu \in \mathfrak J_n} \big\|f_\nu \|_p^p \big)^{1/p}\lesssim\delta_n^{1-\frac{N+1}p}2^{-\frac kp}\|f \|_p.\]
Therefore, taking sum over  $n$,  we get \eqref{overline-a0}, which proves Proposition \ref{newprop}.

It remains to prove \eqref{tannu}.
By Lemma \ref{symcheck}, we have $C^{-1}\mathfrak a_{\nu}^{n} \in \mathfrak A_k(\delta_n,\delta_n\nu)$ for a constant $C>0$.  For $n=0$, it  is easy to show  \eqref{tannu}. Since  $\delta_0=2^{-k/N}$,  applying  Lemma \ref{kernel}, we get
\[ 
\| \mathcal T[\mathfrak a_{\nu}^0]f\|_p\lesssim \delta_0^{1-\frac 1p}\|f\|_p = \delta_0^{1-\frac{N+1}p}2^{-\frac kp}\|f\|_p, \quad 2\le p\le \infty. 
\]
For $n\ge1$, we need to decompose $\mathfrak a_{\nu}^{n}$ further.
Let us set
\[
\mathfrak a_{\nu,1}^{n} (s,t,\tau,\xi) = 
\mathfrak a_{\nu}^{n} (s,t,\tau,\xi) (1-
\beta_0 )\big( 10 \delta_n^{-2N!} |\langle 2^{-k}G(s),(\tau,\xi)\rangle |^{2(N-1)!}\big)
\]
and $\mathfrak a_{\nu,0}^{n} = \mathfrak a_{\nu}^{n}-\mathfrak a_{\nu,1}^{n}$, so  we have
$\mathfrak a_{\nu}^{n} = \mathfrak a_{\nu,1}^{n}+\mathfrak a_{\nu,0}^{n}. $  
We note that  $C^{-1} \mathfrak a_{\nu,i}^{n} \in \mathfrak A_k(\delta_n,\delta_n\nu)$, $i=0,1$ for some $C>0$. This can be shown by following the proof of Lemma \ref{symcheck}. So, we omit the detail. 

We now decompose $\mathcal T [\mathfrak a_{\nu}^{n}]f=\mathcal T [\mathfrak a_{\nu,1}^{n}]f+\mathcal T [\mathfrak a_{\nu,0}^{n}]f.$
For  \eqref{tannu}, it suffices to show
\begin{align}\label{goal}
\|\mathcal T [\mathfrak a_{\nu,i}^{n}]f \|_p \le C \delta_n^{1-\frac {N+1}p} 2^{-\frac kp}\|f \|_p, \quad i=0,1,
\end{align}
for $p>2N-2$. 
It is clear that \eqref{lowscale} holds with $\delta=\delta_n$, $\sz=\delta_n\nu$, and some large $B$ on $\supp \mathfrak a_{\nu,1}^{n}$.  
 By Corollary \ref{lessdegen} we have  \eqref{goal} for $i=1$  if $p>2N-2$. The operator $|\mathcal T [\mathfrak a_{\nu,0}^{n}]$ can be handled in the same manner as
$ \mathcal T[\mathfrak a_1]$ since 
\Be
\label{lll}
|\tau+\langle\gamma(s),\xi\rangle|\gtrsim  \delta_n^{N} 2^{k} 
\Ee
holds on  $ \supp \mathfrak a_{\nu,2}^{n}$. We set $\mathfrak a=-i2^{k}\delta_n^{N}(\tau+\langle \gamma(s), \xi \rangle )^{-1} \partial_{t}\mathfrak a_{\nu,0}^{n}$. 
Integration by parts in $t'$ and \eqref{T} yields $\mathcal T [\mathfrak a_{\nu,0}^{n}]=(2^{k}\delta_n^{N})^{-1}\mathcal T [\mathfrak a]$. Using \eqref{lll} and 
the fact that  $C^{-1} \mathfrak a_{\nu,2}^{n} \in \mathfrak A_k(\delta_n,\delta_n\nu)$  for some $C>0$,  one can easily verify that \eqref{symsupp} and \eqref{symineq2}  hold for $\mathfrak a$ with $\delta=\delta_n$, $\sz=\delta_n\nu$. 
Thus, by Lemma \ref{kernel}, we have
\[
\|\mathcal T [\mathfrak a_{\nu,0}^{n}]f \|_p \lesssim \delta_n^{1-\frac1p} (\delta_n^N2^{k})^{-1}\|f\|_p \lesssim  \delta_n^{1-\frac {N+1}p} 2^{-\frac kp}\|f \|_p
\]
for $p\ge2$, which gives  \eqref{goal} for $i=0$. For the second inequality we use the fact that $\delta_n\ge 2^{-k/N}$.
\qed

\section{Proof of  Theorem \ref{thm:improving}}\label{order-up}

We first prove the sufficiency part, that is to say,  the estimate \eqref{A1} with $\alpha=1-1/p$ for $1\le p<p_d$  by making use of  Theorem \ref{lclb}.

\subsection{Proof of the estimate \eqref{A1} with $\alpha=1-1/p$}\label{sec:FS}
We make use of  the argument in \cite{PS06, PRS}.
As mentioned before, it suffices to prove \eqref{dual2}  by duality.
Let $P_k$ denote the (Littlewood-Paley projection) operator defined by
\[
\mathcal F(P_kg)(\xi,\tau)=\beta(2^{-k}|(\xi,\tau)|)\widehat g(\xi,\tau), \quad k\ge1
\] 
for $\beta \in \mathrm C_0^\infty([1/2,2])$.
Recall that $\beta_0\in \mathrm C_0^\infty([-1,1])$ such that $\beta_0= 1$ on $[-1/2,1/2]$ and set 
$ \beta_\ast(t)=\beta_0(C_0^{-1}2^{-6}t)-\beta_0(C_0 2^{6}t)$.  Here $C_0=1+2\sup \{|\gamma(s)|+|\gamma'(s)|: s \in \supp \psi \}$.
Let  $f_k$ be given by 
\[
\widehat{f_k}(\xi,u)=\beta_\ast(2^{-k} |(\xi,u)|) \widehat f(\xi,u).
\]
We claim that 
\begin{align}\label{L-P}
\big\| \big(\sum_{k \ge1}| P_k \mathcal R f|^2 \big)^{1/2} \big\|_p
\lesssim \big\| \big( \sum_{k\ge1} 2^{-\frac {2k}p} | f_k|^2 \big)^{1/2} \big\|_p+\|f\|_{L_{- M}^p}
\end{align}
for $p>2d$ and $M\gg1$.  Then \eqref{dual2} follows by the Littlewood-Paley inequality.

Let $\tilde \beta=\beta_0(2^{-3}\,\cdot)-\beta_0(C_02^{3}\,\cdot)$
Considering an operator $\mathcal R_k$ given  by
\[
\mathcal F_x (\mathcal R_kf)(\xi,t)=\tilde\beta({|\xi|}/{2^k}) \mathcal F_x(\mathcal Rf)(\xi,t),
\]
we decompose 
\begin{align}\label{express}
P_k \mathcal Rf=P_k \mathcal R_k  f_k+P_k \mathcal R_k (f-f_k)+P_k(\mathcal R-\mathcal R_k) f.
\end{align}
In what follows we show that the contributions from the second and third terms are negligible.  In fact,  for any $M\ge1$ if $p\ge1$, 
we have
\begin{align}\label{Ek}
\big\| \big(\sum_k |P_k \mathcal R_k(f-f_k)|^2 \big)^{1/2}\big\|_p \lesssim \|f\|_{L_{-M}^p}
\end{align}
and \eqref{error-2} below.

 To see \eqref{Ek}, 
%
note
$
\mathcal F_x(\mathcal R_k g)(\xi,t')= 
\int m(\xi,t',u) \widehat g(\xi,u)\,du
$
where 
\[m(\xi,t',u)=(2\pi)^{-1}\chi(t') \tilde\beta({|\xi|}/{2^k})
\int e^{i(s u-t'\gamma(s)\cdot \xi)}\psi(s)\,ds.
\]
Since $|(\xi,u)|\ge C_0 2^{k+5}$ or $|(\xi,u)| \le C_0^{-1}2^{k-5}$ on  $\supp \mathcal F (f-f_k)$, we have $|u| \ge C_0 |\xi|$ if $C_0^{-1}2^{k-4}\le|\xi| \le 2^{k+3}$.
Therefore,   integration by parts gives
\[ |\partial_{\xi,u}^\alpha m(\xi,t',u)| \lesssim 2^{-kN} (1+|(\xi,u)|)^{-N}, \quad (\xi,u)\in  \supp \mathcal F (f-f_k)\] 
for any $\alpha$ and $N\ge1$.
Note that
$
\textstyle P_k \mathcal R_kg(x,t)=\int K(x,y,t,s') g(y,s')\,dyds'
$
where
\[
K(x,y,t,s')=\frac1{(2\pi)^{d+1}}\int  e^{i(x-y,t-t')\cdot(\xi,\tau)} e^{-is'u}\beta\Big(\frac {|(\xi,\tau)|}{2^{k}}\Big) m(\xi,t',u)
\,d\xi d\tau dudt'.
\]
Thus if $g=P_j(f-f_k)$, then $|\partial_{\xi,u}^\alpha m(\xi,t',u)| \lesssim 2^{-kN} 2^{-jN}$ for $(\xi,u) \in \supp \mathcal Fg$ and integration by parts shows
\[|K(x,y,t,s')| \lesssim 2^{-kN}2^{-jN}(1+|x-y|+|s'|)^{-N}(1+|t|)^{-N}.\]
Decomposing  $\mathcal R_k(f-f_k)=\sum_{j}\mathcal R_kP_j(f-f_k)$, 
we get \eqref{Ek} for any $M\ge1$ and $p\ge1$.

We now show
\begin{align}\label{error-2}
\big\| \big( \sum_k | P_k(\mathcal R-\mathcal R_k)f |^2 \big)^{1/2} \big\|_p
\lesssim \|f\|_{L_{-M}^p}
\end{align}
 for $p\ge1$ and $M\ge1$.
We write $\mathcal F(\mathcal Rf-\mathcal R_kf)(\xi,\tau)=\int b(s,\xi,\tau)\mathcal F_xf(\xi,s)\,ds$ where
\[
b(s,\xi,\tau)=\frac1{2\pi}\int e^{it'(\gamma(s)\cdot \xi-\tau)} \big( 1-\tilde\beta({|\xi|}/{2^k})\big) \chi(t')\,dt'  \psi(s).
\]
Since $|\xi| \le C_0^{-1}2^{k-2}$ or $|\xi|\ge 2^{k+2}$ on $\supp \mathcal F_x(\mathcal Rf-\mathcal R_kf)$,
we have $|\tau| \ge C_0 |\xi|$ if $2^{k-1} \le |(\xi,\tau)| \le 2^{k+1}$.  Integration by parts gives $|\partial_\xi^\alpha b(s,\xi,\tau)|\lesssim 2^{-kN}$ for any $\alpha$ and $N$. Hence,  
\begin{align}\label{est01}
\|P_k(\mathcal R -\mathcal R_k) f\|_{p} \lesssim 2^{-kN} \|f\|_p, \quad p\ge 1
\end{align}
 for  all $N\ge1$. 
Since $|\xi| \le C_0^{-1}2^{k-2}$ on $\supp\mathcal F(P_k(\mathcal R-\mathcal R_k)f)$, similarly as in the proof of \eqref{Ek}, we have
$\|P_k(\mathcal R -\mathcal R_k) P_jf\|_{p} \lesssim 2^{-jN} \|P_j f\|_p$ for $j\ge k+C'$ for some $C'\ge1$. The estimate  \eqref{est01} gives $\|P_k(\mathcal R -\mathcal R_k) P_jf\|_{p} \lesssim 2^{-kN} \|P_j f\|_p$ for $j\le k+C'$. Combining those estimates, we get \eqref{error-2}.

Therefore,    the estimate \eqref{L-P} follows if we  show 
\begin{align}\label{square}
\big\| \big(\sum_k|P_k \mathcal R_k f_k|^2 \big)^{1/2} \big\|_p
\lesssim \big\| \big( \sum_{k\ge1} 2^{-\frac {2k}p} | f_k|^2 \big)^{1/2} \big\|_p
\end{align}
for $p>2d$.
This can be done by using \cite[Theorem 1]{PRS} and \eqref{smoothing} (also see \cite{PS06, PS07,BGHS}).  Indeed, let $ \tilde \beta \in \mathrm C_c^\infty((1/4,4))$ such that 
$\tilde \beta\beta=\beta$. Consider the operator $\tilde P_k$ given  by 
$\mathcal F(\tilde P_kg)(\xi,\tau)=\tilde \beta(2^{-k}|(\xi,\tau)|) \widehat g(\xi,\tau)$.  
Note that $P_k\mathcal R_k f_k=P_k \tilde P_k \mathcal R_k f_k$.

Let us denote  the center of  a cube $Q$ by $(x_Q,t_Q)$ and
 set \[\mathcal E_Q=\{(y,s): \dist( y-x_Q, t_Q \gamma(I)) \le 10 \text{diam} (Q), ~s \in I\}.\]
 Since $T_k=\tilde P_k\mathcal R_k$ and $\mathcal E_Q$ satisfy the assumptions in \cite[Theorem 1]{PRS}, by using \eqref{smoothing} we obtain \eqref{square}.
We omit the details.
\qed

\subsection{Sharpness of  smoothing order } 
In this section, we show upper bounds on the smoothing order $\alpha$ for which $L^p$--$L_\alpha^p$ estimate for $\mathfrak Rf$ holds when $\gamma$ is of maximal type $L$. 
In  \cite{PS06} those  bounds were obtained for $d=2$.   Modifying the examples in  \cite{PS06}, 
we show the following.  

\begin{prop}\label{nece} 
Let $d\ge3$, $L\ge d$, and $1 \le p \le \infty$. Let $\psi$ and $\chi$ be nontrivial, nonnegative continuous functions supported in the interiors of  $I$ and $[1,2]$, respectively. 
Suppose there is an $\sz$ such that $\psi(\sz)\neq0$ and $\gamma$ is of type $L$ at $\sz$. 
Then, $\mathfrak Rf$ maps $L^p(\mathbb R^{d+1})$ boundedly to $L_\alpha^p(\mathbb R^{d+1})$  only if 
\[ (i)\ \alpha\le  1-p^{-1},  \quad (ii)\    \alpha\le  (2d)^{-1}, \quad (iii)\    \alpha\le (Lp)^{-1} .\]
\end{prop}

In particular, the upper bound $(i)$ provides 
the necessity part of Theorem \ref{thm:improving}, thus, 
the proof Theorem \ref{thm:improving} is completed.  We prove  the upper bounds $(i)$, $(ii)$, and $(iii)$, separately. 


\begin{proof}[Proof of $($i\,$)$] 
Let $t_{\mathsmaller 0} \in (1,2)$ such that $\chi(t_{\mathsmaller 0})>0$.  
We choose $\zeta \in \mathcal S(\mathbb R^{d})$ such that $\zeta \ge1$ on $[-1,1]^{d}$,
$\supp \widehat \zeta \subset [1/2,4]^d$, and $\widehat\zeta=1$ on $[1,2]^d$.
Let $\psi_0\in \mathrm C_c^\infty((-1,1))$ satisfy $\psi_0=1$ on $[-1/2,1/2]$.
We take
\[ 
f(x,t)=\zeta(\lambda x)\psi_0(\lambda r_0 |t-t_{\mathsmaller 0}|), 
\] 
where $r_0=1+\sup_{s \in I}|\gamma(s)|$.
Note $\mathfrak Rf(x,s) \gtrsim \lambda^{-1}$  if  $|x+t_{\mathsmaller 0}\gamma(s)|\le c \lambda^{-1}$ and $|s-\sz|<c$ for a small constant $c>0$.
Thus,  $\| \mathfrak Rf \|_{L^p(\mathbb R^{d+1})} \gtrsim \lambda^{-1-d/p}$.
Since  
\[
\mathcal F_x(\mathfrak Rf(\cdot, s))(\xi) = \lambda^{-d} \psi(s)\int \widehat \zeta(\lambda^{-1}\xi) e^{it\gamma(s)\cdot \xi}\psi_0(\lambda|t-t_{\mathsmaller 0}|)\chi(t)\,dt, 
\] 
it follows that $\supp_\xi \mathcal F_x(\mathfrak Rf)$ is included in $\{ \xi : |\xi| \sim \lambda \}$. Hence, $\| \mathfrak Rf(\cdot,s) \|_{L_\alpha^p(\mathbb R^{d};dx)} \gtrsim \lambda^{\alpha-1-d/p}$, so 
we have $\| \mathfrak Rf\|_{L_\alpha^p(\mathbb R^{d+1})} \gtrsim \lambda^{\alpha-1-d/p}$.
Since $\|f\|_p \lesssim \lambda^{-(d+1)/p}$, we get $\alpha \le 1-1/p$.
\end{proof}

\begin{proof}[Proof of $($ii\,$)$] 
Let  $\tilde I \subset (-1,1)$ be a nonempty compact  interval  such that \eqref{nonv} holds for $s \in \tilde I$.
Also, we fix a constant $\rho \gg1$ to be chosen later.
Let $\{s_\ell\}\subset \tilde I$ be a collection of $\rho\lambda^{-1/d}$-separated   points which are as many as $C\rho^{-1}\lambda^{1/d}$.  
Since $G(s_\ell), G'(s_\ell), \dots, G^{(d-1)}(s_\ell)$ are linearly independent in $\mathbb R^{d+1}$,  there is a  unit vector $\Xi_\ell \in \big(  \sspan\{ G^{(j)}(s_\ell): j=0,1,\dots,d-1\} \big)^\perp$.

Let $\phi \in \mathcal S(\mathbb R^{d+1})$ such that $\phi \ge1$ on $[-3r_0,3r_0]^{d+1}$
 and $\widehat \phi$ is supported in $[-1,1]^{d+1}$ where $r_0=1+\sup_{s \in I} |\gamma(s)|$. 
Let $\varepsilon_\ell \in \{\pm 1\}$ be independent random variables. We consider
\[
\textstyle f(x,t)=\sum_\ell \varepsilon_\ell  f_\ell(x,t):=   \sum_\ell \varepsilon_\ell \phi(x,t) e^{i\lambda \Xi_\ell \cdot (t,x)}.
\]

Since $\langle \Xi_\ell, G^{(j)}(s_\ell) \rangle=0$ for $j=0,\dots,d-1$,
by Taylor's theorem we have 
\begin{align}\label{est2}
\langle \Xi_\ell, G(s) \rangle =\langle \Xi_\ell, G^{(d)}(s_\ell) \rangle(s-s_\ell)^{d}/d!+O(|s-s_\ell|^{d+1}).
\end{align}
Thus $|t \langle \Xi_\ell , G(s)\rangle | \le 2^{-2}\lambda^{-1}$ whenever 
$s \in I_\ell:=\{ s \in \tilde I : |s-s_\ell| \le c \lambda^{-1/d}\}$ 
for a  $c>0$ small enough.
Noting that
\begin{align}\label{est1}
\mathfrak Rf_\ell(x,s)= e^{i\lambda \Xi_\ell \cdot (0,x)} \psi(s)\int \phi (x+t\gamma(s),t ) e^{i\lambda t\Xi_\ell \cdot G(s)} \chi(t)\,dt,
\end{align}
we see 
$|\mathfrak Rf_\ell(x,s)| \gtrsim1$ if $ (x,s) \in B_\ell:=[-c,c]^d\times I_\ell.$  Thus, 
$
\sum_\ell \big\| \mathfrak Rf_\ell \big\|_{L^p(B_\ell)}^p \gtrsim \rho^{-1}.
$
Meanwhile, by \eqref{est1}, \eqref{est2},  and  integration by parts in $t$  we have
$
|\mathfrak Rf_m(x,s)| \lesssim (1+\lambda |s_\ell-s_m|^{d})^{-N}
$
for any $N\ge1$ if $m\neq \ell$ and $s \in I_\ell$. 
Since $\{s_\ell\}$ are $\rho \lambda^{-\frac{1}{d}}$-separated, it is easy to see
\[
\sum_\ell \big\| \sum_{m\neq \ell} |\mathfrak Rf_m| \big\|_{L^p(B_\ell)}^p \lesssim \sum_\ell \sum_{m \neq \ell} (1+\lambda |s_\ell-s_m|^{d})^{-pN}  \lambda^{- 1/d} \lesssim \rho^{-pdN-1}.
\]
Therefore, taking $\rho,N$ sufficiently large,   we have $\|\mathfrak Rf \|_p \gtrsim \rho^{-1}$ for any choice of $\varepsilon_\ell$.

By our choice of $\phi$ it follows that  $\mathcal F_x (\mathfrak Rf)$ is supported on $\{\xi : C_1\lambda \le |\xi| \le  C_2 \lambda \}$ for some 
positive constant $C_1, C_2$. Thus, 
$\| \mathfrak Rf \|_{L_\alpha^p} \gtrsim  \lambda^\alpha\|\mathfrak Rf\|_p $. Combining this with the $L^p$--$L_\alpha^p$ estimate  
gives 
$\lambda^\alpha   \le C  \|f\|_p $
for any choice of $\varepsilon_\ell$. By Khinchine's inequality we have  $
\mathbb E (\| f \|_p^p) \sim \int (\sum_\ell | f_\ell|^2 )^{\frac p2}\,dxdt \sim C_\rho \lambda^{\frac{p}{2d}}.$  
Therefore, we see $\lambda^\alpha\lesssim \lambda^\frac1{2d}$ and then $\alpha \le 1/(2d)$ taking $\lambda\to \infty$.
\end{proof}

\begin{proof}[Proof of $($iii\,$)$] 
Since $\gamma$ is of type $L$ at $\sz$, by an affine transformation and taking $\psi$ supported near $\sz$, we may assume 
\[
\gamma(s+\sz)=\gamma(\sz)+(s^{a_1} \varphi_1(s),\dots, s^{a_d} \varphi_d(s))
\]
for $1\le a_1 <\dots<a_d=L$  and  smooth functions $\varphi_j$, $j=1,\dots,d$, 
where  $\|\varphi_j-1/a_j!\|_{\mathrm C^{a_d+1}(I)}\le c$ for a small constant $c>0$. 
We may also assume $\sz=0$ and furthermore $\gamma(0)=0$ by replacing $f(x,t)$ by $f(x-t\gamma(0),t)$.

Let $\phi_1\in \mathcal S(\mathbb R)$ such that $\phi_1 \ge1$ on $[-1,1]$,
and $\supp \widehat \phi_1\subset [1/2,4]$ with $\widehat\phi_1=1$ on $[1,2]$. 
Let $\psi_0\in \mathrm C_c^\infty((-1,1))$ with $\psi_0=1$ on $[-1/2,1/2]$. We consider 
\[
f(x,t)=\prod_{j=1}^{d-1} \psi_0 ( \lambda^{a_j/L}x_j )
\phi_1 ( \lambda x_{d} ) \chi(t).
\]
Denoting $\|a\|=\sum_{j=1}^d a_j$, we have $\|f\|_p \lesssim \lambda^{-\|a\|/(Lp)}$.
Set
$
E_\lambda=\big\{ (x,s) \in \mathbb R^d \times I: |x_j| \le c \lambda^{-a_j/L}, j=1,\dots,d, \quad~ |s| \le c \lambda^{-1/L} \big\}
$
for a sufficiently small $c>0$. 
Since $\gamma(s)=(s^{a_1} \varphi_1(s),\dots, s^{a_d} \varphi_d(s))$, 
$
| \langle x+t\gamma(s),e_j \rangle| \le 2^{-1}\lambda^{-a_j/L},$ $j=1,\dots,d$, 
for  $(x,s) \in E_\lambda$ and $t\in [1,2]$. So,   $\mathfrak Rf (x,s)\gtrsim 1$ for  $(x,s) \in E_\lambda$. This gives  
$\| \mathfrak Rf\|_p\gtrsim \lambda^{-(\|a\|+1)/(Lp)}$.
Since $\supp \mathcal F_{x_{d}} (\mathfrak Rf) \subset \{ \xi_{d} : |\xi_{d}| \sim \lambda\}$, 
$
\| \mathfrak Rf \|_{L_\alpha^p} \gtrsim \lambda^{\alpha-(\|a\|+1)/(Lp)}.
$
Therefore,  we obtain  $\alpha \le 1/(Lp)$.
\end{proof}

\subsection*{Acknowledgement} The research  was  supported by NRF2022R1A4A1018904 (Ko, Lee, and Oh) and  NRF2022R1I1A1A01055527 (Ko), and a KIAS Individual Grant MG089101 (Oh).  

\end{document}